\documentstyle[12pt,a4,pstricks,pst-node,pst-coil]{article}

\catcode`@=11
\@addtoreset{equation}{section}
\catcode`@=12

\newcommand{\nc}{\newcommand}

\font\Bbbtwl=msbm10 scaled 1200
\newfam\Bbbfam
\textfont\Bbbfam=\Bbbtwl
\font\eufmtwl=eufm10 scaled 1200
\newfam\eufmfam
\textfont\eufmfam=\eufmtwl

\newfont{\Bbb}{msbm10 scaled 1200}
\newfont{\goth}{eufm10 scaled 1200}
\newfont{\typewriter}{cmtt10 scaled 1200}

\nc{\Spinc}{\mbox{{$\mbox{spin}^c$}}} %neu!
\nc{\Spinh}{\mbox{{$\mbox{spin}^h$}}} %neu!
\nc{\specDi}{\mbox{{$\mbox{spec}\,(D_i)$}}}%neu!
\nc{\specDe}{\mbox{{$\mbox{spec}\,(D_1)$}}}%bernommen
\nc{\specDz}{\mbox{{$\mbox{spec}\,(D_2)$}}}%bernommen
\newcommand{\RR}{\mbox{\Bbb R}}%bernommen
\newcommand{\CC}{\mbox{\Bbb C}}%bernommen
\newcommand{\ZZ}{\mbox{\Bbb Z}}%bernommen

\nc{\NN}{\mbox{{${\rm I\!N}$}}}
\nc{\zz}{\mbox{{$\scriptstyle{\rm Z\!\!Z}$}}}
\nc{\nn}{\mbox{$\scriptstyle{\rm I\!N}$}}
\nc{\adick}{\mbox{{$\bar{a}$}}}
\nc{\bdick}{\mbox{{$\bar{b}$}}}
\nc{\Zdick}{\mbox{{$\bar{Z}$}}}
\nc{\Alphadick}{\mbox{{$\bar{\alpha}$}}}
\nc{\Cl}{\mbox{{${\it Cl}$}}}
\nc{\SOn}{\mbox{{$\mbox{SO}(n)$}}}
\nc{\SOzwoN}{\mbox{{$\mbox{SO}(2N)$}}}
\nc{\sozwom}{\mbox{{$\mbox{so}(2m)$}}}
\nc{\Spinn}{\mbox{{$\mbox{Spin}(n)$}}}
\nc{\SpinzwoN}{\mbox{{$\mbox{Spin}(2N)$}}}
\nc{\GLplusn}{\mbox{{$\mbox{GL}^+(n)$}}}
\nc{\PsoM}{\mbox{{$P_{\mbox{\scriptsize SO}}\;M$}}}
\nc{\PspinM}{\mbox{{$P_{\mbox{\scriptsize spin}}\;M$}}}
\nc{\PspinMN}{\mbox{{$P_{\mbox{\scriptsize spin}}\;M|_N$}}}
\nc{\dhX}{\mbox{{$(\mbox{\it dh} \cdot X)$}}}
\nc{\dotrho}{\mbox{{$\dot{\rho}$}}}
\nc{\dotgamma}{\mbox{{$\dot{\gamma}$}}}
\nc{\doppeldotgamma}{\mbox{{$\ddot{\gamma}$}}}
\nc{\dotB}{\mbox{{$\dot{B}$}}}
\nc{\s}{\mbox{{$\mbox{\bf s}$}}}
\nc{\ce}{\mbox{{$\mbox{\bf c}$}}}
\nc{\sint}{\mbox{{$\mbox{sin}(t)$}}}
\nc{\cost}{\mbox{{$\mbox{cos}(t)$}}}
\nc{\sinht}{\mbox{{$\mbox{sinh}(t)$}}}
\nc{\sinheins}{\mbox{{$\mbox{sinh}(1)^2$}}}
\nc{\cosht}{\mbox{{$\mbox{cosh}(t)$}}}
\nc{\cosheins}{\mbox{{$\mbox{cosh}(1)^2$}}}
\nc{\sinhC}{\mbox{{$\mbox{sinh}(C)$}}}
\nc{\coshC}{\mbox{{$\mbox{cosh}(C)$}}}
\nc{\sinhlt}{\mbox{{$\mbox{sinh}(|\lambda| t)$}}}
\nc{\coshlt}{\mbox{{$\mbox{cosh}(|\lambda| t)$}}}
\nc{\dt}{\mbox{\it dt}}
\nc{\specD}{\mbox{{$\mbox{spec}\,(D)$}}}
\nc{\specDT}{\mbox{{$\mbox{spec}\,(D_T)$}}}
\nc{\spectildeDT}{\mbox{{$\mbox{spec}\,(\tilde{D}_T)$}}}
\nc{\specDstrich}{\mbox{{$\mbox{spec}\,(D')$}}}
\nc{\Hom}{\mbox{{$\mbox{Hom}\,$}}}
\nc{\detDgamma}{\mbox{{$\mbox{det}\, D_{\gamma}$}}}
\nc{\Fraglichvee}{\mbox{{$V$}}}
\nc{\overRAunder}{\mbox{${\mathop{\longrightarrow}\limits^{\theta_{*}}_{\cong}}$}}
\nc{\pigstern}{(\pi_{\gamma})_*}
\nc{\vsp}{\vspace}
\nc{\RE}{\mbox{\bf \goth Re}}
\nc{\dm}{\mbox{dim}}

\begin{document}

\psset{unit=0.4cm}
\title{Harmonic Spinors for Twisted Dirac Operators}
\author{Christian B\"ar\thanks{partially supported by SFB 256 and by the 
GADGET program of the EU}}
\date{May, 1996}
\maketitle

\begin{abstract}
\noindent
We show that for a suitable class of ``Dirac-like'' operators there holds 
a Gluing Theorem for connected sums. More precisely, if $M_1$ and $M_2$
are closed 
Riemannian manifolds of dimension $n\ge 3$ together with such operators, 
then the connected sum $M_1 \# M_2$ can be given a Riemannian metric such 
that the spectrum of its associated operator is close to the disjoint 
union of the spectra of the two original operators. As an application, we 
show that in dimension $n\equiv 3$ mod 4 harmonic spinors for the Dirac 
operator of a spin, $\Spinc$, or $\Spinh$ manifold are not 
topologically obstructed. 

\end{abstract}

\section{Introduction}

The interplay between geometric, topological, and analytic invariants of
Riemannian manifolds is one of the major topics of current investigations in 
global ana\-lysis and differential geometry.
A classical example is provided by Hodge-deRham theory.
The dimension of the kernel of the Laplace-Beltrami operator acting on 
$p$-forms on a closed Riemannian manifold is an important topological
invariant, the $p$-th Betti number.
In particular, it does not depend on the Riemannian metric.
The question arises whether one can obtain further topological
invariants using other natural
operators like the Dirac operator on a closed Riemannian spin manifold.

There are topological restrictions against existence of harmonic spinors
in dimension 2.
A 2-sphere equipped with an arbitrary Riemannian metric $g$ does not
have non-trivial harmonic spinors.
This can be deduced e.g.\ from the eigenvalue estimate \cite[Thm.\ 2]{b3}
$$
\lambda^2 \geq \frac{4\pi}{\mbox{area}(S^2,g)}
$$
which holds for all Dirac eigenvalues $\lambda$ on $(S^2,g)$.
At the moment the 2-sphere is the only closed spin manifold for which
one knows non-existence of harmonic spinors for {\em all} Riemannian
metrics and all spin structures.
In fact, the present paper gives some evidence to the conjecture that
$S^2$ is the only such manifold.

On a surface of genus 1 or 2 the dimension of the space of harmonic
spinors does not depend on the Riemannian metric, but it does depend
on the choice of spin structure.
For genus larger than 2 it depends on both, the Riemannian metric and
the spin structure.
One can always choose the metric and the spin structure in such a way
that there are non-trivial harmonic spinors \cite{hi}.
For hyperelliptic surfaces one can compute the dimension for all
spin structures \cite{bs}.

Hitchin \cite{hi} computed the Dirac spectrum for a suitable
1-parameter family of metrics on $S^3$. 
It turns out that for generic parameter values there are no
non-trivial harmonic spinors but for special choices of the
parameter the dimension of the space of harmnonic spinors becomes
arbitrarily large.

If the dimension of the manifold is divisible by 4, then the
Atiyah-Singer Index Theorem \cite[Thm.\ 5.3]{asIII} implies
$$
\dim\{\mbox{harmonic spinors}\} \geq |\hat{A}(M)|
$$
where $\hat{A}(M)$ is a topological invariant, the $\hat{A}$-genus of $M$.
Kotschick \cite{ko} exhibited algebraic surfaces for which 
$\dim\{\mbox{harmonic spinors}\}$ exceeds $|\hat{A}(M)|$ arbitrarily much.

All known examples indicate that the following conjecture should hold.

\vspace{1cm}
{\bf Conjecture.}
{\em Let $M$ be a closed spin manifold of dimension $n \geq 3$.
Let the spin structure on $M$ be fixed.

Then there exists a Riemannian metric on $M$ such that there are
non-trivial harmonic spinors.}

\vspace{1cm}
Hitchin \cite{hi} proved this conjecture by differential 
topological methods in dimension $n \equiv 0,\pm 1$ mod 8.
In \cite{b4} we showed by analytic methods that the conjecture also
holds in dimension $n \equiv 3$ mod 4.

In this paper which should be regarded as a sequel to \cite{b4} we enlarge
the class of operators for which the analogous statement holds.
This is possible because the methods of \cite{b4} are essentially local.
This class of operators contains in particular the Dirac operators of 
$\Spinc$ or $\Spinh$ manifolds which have regained interest recently 
because of Seiberg-Witten theory \cite{w}, see \cite{morgan,salamon} 
for an introduction.
Some care has to be taken with the statement of the result.
In contrast to spin manifolds the Dirac operator of a $\Spinc$ manifold
is not determined by the Riemannian metric alone but also depends on the
choice of a connection on the canonical line bundle.
A similar remark holds for $\Spinh$ manifolds.
The nature of the Seiberg-Witten equations indicates that one
should look at modifications of the metric and of the connection on the
canonical line bundle separately.
It turns out that a modification of the metric is 
sufficient to produce a non-trivial kernel.
Theorem 4.1 applied to $\Spinc$ manifolds says

\vspace{1cm}
{\bf Theorem.}
{\em Let $M$ be a closed $\Spinc$ manifold of dimension $n \equiv$ 3 mod 4.
Let a connection on the canonical line bundle be fixed.
Then there exists a Riemannian metric on M such that there are non-trivial
harmonic spinors for the associated Dirac operator.}

\vspace{1cm}
\noindent
For $n=3$ this can be rephrased in terms of Seiberg-Witten equations
for 3-dimensional (!) manifolds.
Of course, it would be interesting to have the analogous statement in
dimension 4, but at the moment this seems out of reach.

The construction shows that this metric can be obtained by deforming
any given metric in an arbitrarily small open set while keeping it 
unchanged outside.
Since the construction is local Theorem 4.1 applies to
self-adjopint elliptic operators which look like a twisted Dirac operator 
in {\em some} non-empty open subset of the manifold.
Outside this set the operator can be anything and will not be modified.
The precise statement is as follows

\vspace{1cm}
{\bf Theorem 4.1.} {\it Let} $M$ {\it be a closed Riemannian manifold of 
dimension $n\equiv$ 3 mod 4.} 
{\it Let $D$ be an elliptic self-adjoint differential operator over $M$ 
of order 1.} 
{\it Let $U\subset M$ be a non-empty open subset}. 
{\it Let the restriction of $D$ to $U$ be a twisted Dirac operator}.

{\it Then one can deform the Riemannian metric in $U$ such that the 
resulting operator $\tilde{D}$ has non-trivial kernel}.

\vspace{1cm}
\noindent                                        
The proof relies firstly on the computation of the spectrum of the classical 
Dirac operator for a certain 1-parameter family of metrics on odd-dimensional
spheres, the so-called {\em Berger metrics}.
This computation has been carried out in \cite{b4}.

Secondly, we prove a gluing theorem for such operators on connected sums.
Given two closed manifolds of dimension $n \geq 3$ together with
operators which look like the classical Dirac operator (or a multiple
of it) in some non-empty open subset, then by removing balls in these
subsets and gluing one can form the connected sum with a Riemannian metric 
such that the spectrum of the associated operator on the sum is close to 
the disjoint union of the spectra of the two original operators, at
least in some bounded range.

\vspace{1cm}
\newpage

{\bf Theorem 2.1.} ({\it Gluing Theorem})

\noindent
{\it Let $M_1$ and $M_2$ be $n$-dimensional closed 
Riemannian manifolds of dimension $n\ge 3$. Let $U_i \subset M_i$ be 
open balls, let $D_i$ be self-adjoint elliptic differential operators of 
first order over $M_i$ of Dirac type over $U_i$. Let $\Lambda >0$ such 
that $\pm \Lambda \notin \specDe \cup \specDz$. Let $\epsilon >0$. 

Then there exists a Riemannian metric on $X=M_1\# M_2$ and a self-adjoint 
elliptic first order differential operator $D$ over $X$ such that $X$ is 
a disjoint union $X=X_1 \dot\cup X_2 \dot\cup X_3$ where}
\begin{itemize}
\item[(i)]
$X_1$ {\it is isometric to} $M_1 - U_1$ {\it and} $D$ {\it coincides 
with} $D_1$ {\it over} $X_1$,
\item[(ii)]
$X_2$ {\it is isometric to} $M_2 - U_2$ {\it and} $D$ {\it coincides 
with} $D_2$ {\it over} $X_2$,
\item[(iii)]
$X_3$ {\it is diffeomorphic to} $(0, 1) \times S^{n-1}$ {\it and} $D$ {\it 
is of Dirac type over} $X_3$, 
\end{itemize}
{\it and such that }
$D$ {\it is} $(\Lambda, \epsilon)$-{\it spectral close to the disjoint 
union} $D_1\dot\cup D_2$ {\it over} $M_1\dot\cup M_2$.

\vspace{1cm}
\noindent
In \cite{b4} we gave a proof of this Gluing Theorem for the classical
Dirac operator on odd dimensional spin manifolds, $n\geq 3$.
The restriction to odd dimension had technical reasons.
We worked with explicit solutions of the eigenspinor equation on
Euclidean annuli which can be given by a power series if the
dimension is odd.
In even dimensions additional logarithmic terms appear.
Although it is likely that one can carry over the proof of \cite{b4}
we chose to avoid the use of explicit solutions in this paper
and we work with a-priori estimates instead.

In Section 2 we first describe the class of operators under consideration
and then we formulate the Gluing Theorem.
The proof is carried out except for the a-priori estimates on the 
distribution of the $L^2$-norm of eigenspinors on Euclidean annuli.
They are derived in Section 3. 
In Section 4 we apply the Gluing Theorem and prove existence of
metrics with harmonic sections.

An excellent introduction to the Dirac operator of spin and $\Spinc$ 
manifolds is given in \cite{lm} or in \cite{bgv}.
Spin$^h$ manifolds are explained in \cite{b5,na}.
The variational characterization of eigenvalues is explained in \cite{ch},
at least for the Laplace operator.

{\bf Acknowledgement.}
This paper was finished while I enjoyed the hospitality of SFB 288
at Humboldt-Universit\"at zu Berlin.
I would like to thank the referee for several valuable hints on how to
improve the presentation of the results.

\vspace{2cm}
\section{The Gluing Theorem}
The aim of this section is to formulate, and partially prove, a theorem 
which allows us to compare the spectrum of twisted Dirac operators on two 
closed manifolds $M_1$ and $M_2$ with the spectrum of the corresponding 
operator on the connected sum $M_1 \# M_2$ equipped with a suitable metric.

To start, we describe the class of operators that we will consider. If 
$M$ is a closed Riemannian spin manifold, then there is a natural 
self-adjoint elliptic first order differential operator $D$, the {\it 
Dirac operator}. It acts on complex spinor fields and has discrete 
spectrum. Given an additional complex vector bundle $E$ over $M$ with 
connection, one can also form the twisted Dirac operator $D^E$ acting on 
spinors with coefficients in $E$. See \cite{lm} for details. 
Since most of our 
considerations will be local in nature we will be able to avoid global 
topological conditions on $M$ like the spin condition. We make the following

\vspace{1cm}
{\bf Definition.} Let $M$ be a closed Riemannian manifold, let $U\subset 
M$ be open. 
A (formally) self-adjoint elliptic first order differential operator $D$ 
acting on sections of a complex vector bundle over $M$ will be called 
{\it of Dirac type over} $U$, if $D$ coincides over $U$ with a multiple 
of the Dirac operator acting on spinors, i.e. $D$ coincides over $U$ with 
the Dirac operator twisted by a trivial flat bundle.

\vspace{1cm}
\noindent
If the manifold is spin, we can take the Dirac operator itself as an 
example, possibly somehow deformed outside $U$. Another example is given 
by the Dirac operator of a $\Spinc$ manifold provided the canonical line 
bundle is trivial and flat over $U$. 
To get a convenient formulation of the Gluing 
Theorem we introduce the following terminology.

\vspace{1cm}
{\bf Definition.} Let $D_1$ and $D_2$ be two self-adjoint operators with 
discrete spectrum $\specDi$, let $\Lambda>0$, and $\epsilon>0$. We say 
that $D_1$ and $D_2$ are $(\Lambda, \epsilon)$-{\it spectral close} iff 
\begin{itemize}
\item[(i)]
$+\Lambda$, $-\Lambda \notin \specDe \cup \specDz$
\item[(ii)]
$D_1$ and $D_2$ have the same number of eigenvalues in the interval 
$(-\Lambda, \Lambda)$ (counted with multiplicities).
Write 
\begin{eqnarray*}
\specDe \cap (-\Lambda, \Lambda) &=& 
\left\{\lambda_1\le\lambda_2\le\ldots\le\lambda_k\right\}\\
\specDz \cap (-\Lambda, \Lambda) &=&
\left\{\mu_1\le\mu_2\le\ldots\le\mu_k\right\}
\end{eqnarray*}
\item[(iii)]
$|\mu_j - \lambda_j| < \epsilon$ for $j=1, \ldots, k$.
\end{itemize}
Now we can formulate the main result of this section.

\vspace{1cm}
{\bf Theorem 2.1.} ({\it Gluing Theorem})

\noindent
{\it Let $M_1$ and $M_2$ be $n$-dimensional closed 
Riemannian manifolds of dimension $n\ge 3$. Let $U_i \subset M_i$ be 
open balls, let $D_i$ be self-adjoint elliptic differential operators of 
first order over $M_i$ of Dirac type over $U_i$. Let $\Lambda >0$ such 
that $\pm \Lambda \notin \specDe \cup \specDz$. Let $\epsilon >0$. 

Then there exists a Riemannian metric on $X=M_1\# M_2$ and a self-adjoint 
elliptic first order differential operator $D$ over $X$ such that $X$ is 
a disjoint union $X=X_1 \dot\cup X_2 \dot\cup X_3$ where}
\begin{itemize}
\item[(i)]
$X_1$ {\it is isometric to} $M_1 - U_1$ {\it and} $D$ {\it coincides 
with} $D_1$ {\it over} $X_1$,
\item[(ii)]
$X_2$ {\it is isometric to} $M_2 - U_2$ {\it and} $D$ {\it coincides 
with} $D_2$ {\it over} $X_2$,
\item[(iii)]
$X_3$ {\it is diffeomorphic to} $(0, 1) \times S^{n-1}$ {\it and} $D$ {\it 
is of Dirac type over} $X_3$, 
\end{itemize}
{\it and such that }
$D$ {\it is} $(\Lambda, \epsilon)$-{\it spectral close to the disjoint 
union} $D_1\dot\cup D_2$ {\it over} $M_1\dot\cup M_2$.

\vspace{1cm}
\begin{center}
\begin{pspicture}(-10,-10)(10,10)

\pscustom[fillstyle=solid,fillcolor=gray,linestyle=none]{
  \pscurve(-2,7.2)(-0.7,5.8)(-0.7,4.2)(-2,2.8)
  \pscurve[liftpen=2](-2,7.2)(-2.5,6)(-2.5,4)(-2,2.8)}
\psellipse(-5,5)(4.5,3)
\pscurve(-2,7.2)(-2.5,6)(-2.5,4)(-2,2.8)
\pscurve(-6.5,4.8)(-6.2,4.5)(-4.8,4.5)(-4.5,4.8)
\pscurve(-6.2,4.5)(-5.7,5)(-5.3,5)(-4.8,4.5)
\rput(-7.5,5.5){$M_1$}
\rput(-1.6,5){\psframebox*[framearc=0.5]{$U_1$}}

\pscustom[fillstyle=solid,fillcolor=gray,linestyle=none]{
  \pscurve(-2,-2.8)(-0.7,-4.2)(-0.7,-5.8)(-2,-7.2)
  \pscurve[liftpen=2](-2,-2.8)(-2.5,-4)(-2.5,-6)(-2,-7.2)}
\psellipse(-5,-5)(4.5,3)
\pscurve(-6.5,-5.2)(-6.2,-5.5)(-4.8,-5.5)(-4.5,-5.2)
\pscurve(-6.2,-5.5)(-5.7,-5)(-5.3,-5)(-4.8,-5.5)
\rput(-7.5,-4.5){$X_1$}
\pscurve(-2,-2.8)(-2.5,-4)(-2.5,-6)(-2,-7.2)

\pscustom[fillstyle=solid,fillcolor=gray,linestyle=none]{
  \pscurve(2,7.9)(0.8,6.5)(0.8,3.5)(2,2.1)
  \pscurve[liftpen=2](2,7.9)(3,6)(3,4)(2,2.1)}
\psellipse(5,5)(4.5,4)
\pscurve(2,7.9)(3,6)(3,4)(2,2.1)
\pscurve(4.5,4.8)(4.8,4.5)(6.2,4.5)(6.5,4.8)
\pscurve(4.8,4.5)(5.3,5)(5.7,5)(6.2,4.5)
\rput(7,6){$M_2$}
\rput(1.8,5){\psframebox*[framearc=0.5]{$U_2$}}

\pscustom[fillstyle=solid,fillcolor=gray,linestyle=none]{
  \pscurve(2,-2.1)(0.8,-3.5)(0.8,-6.5)(2,-7.9)
  \pscurve[liftpen=2](2,-2.1)(3,-4)(3,-6)(2,-7.9)}
\psellipse(5,-5)(4.5,4)
\pscurve(2,-2.1)(3,-4)(3,-6)(2,-7.9)
\pscurve(4.5,-5.2)(4.8,-5.5)(6.2,-5.5)(6.5,-5.2)
\pscurve(4.8,-5.5)(5.3,-5)(5.7,-5)(6.2,-5.5)
\rput(7,-4){$X_2$}

\pscustom{
  \psecurve(-1,-3.9)(-0.7,-4.2)(-0.3,-4.5)(0.3,-4.5)(0.9,-3.5)(1.1,-3.2)
  \gsave
    \pscurve[liftpen=1](0.9,-6.5)(0.3,-5.5)(-0.3,-5.5)(-0.7,-5.8)
    \fill[fillstyle=solid,fillcolor=gray]
  \grestore}
\psecurve(1.1,-6.8)(0.9,-6.5)(0.3,-5.5)(-0.3,-5.5)(-0.7,-5.8)(-1,-6.1)
\rput(1.8,-5){\psframebox*[framearc=0.5]{$X_3$}}

\end{pspicture}
\begin{center}
{\bf Fig. 1}
\end{center}
\end{center}

\vspace{1cm}
\noindent
The operators that we have in mind for applications are locally twisted 
Dirac operators which are not necessarily globally twisted Dirac 
operators because the manifold need not be spin. Examples are the Dirac 
operators for $\Spinc$ or $\Spinh$ manifolds.

It is an easy exercise to check that the Gluing Theorem does not
hold in dimension 1, e.g.\ for the Dirac operator $D=i\cdot\frac{d}{dt}$
on $S^1_L = \RR /L\cdot\ZZ$, $L > 0$.

\vspace{1cm}
\noindent
The rest of this and the next section are devoted to the proof of the 
Gluing Theorem. To start, we note that \cite[Prop. 7.1]{b4} and its proof 
show the following:

There exists a neighborhood of the Riemannian metric $g_i$ on $M_i$ in 
the $C^1$-topology such that for every metric $g'_i$ in this 
neighborhood which coincides with $g_i$ outside $U_i$ the corresponding 
operator $D'_i$ (which coincides with $D_i$ outside $U_i$ and is of Dirac 
type over $U_i$ for the metric $g'_i$) and $D_i$ are $\left(\Lambda, 
\frac{\epsilon}{2}\right)$-spectral close, say. 
The Taylor expansion of a Riemannian metric in exponential coordinates 
shows that one can deform $g_i$ in $U_i$ such that $U_i$ contains a small 
$n$-ball with Euclidean metric and this deformation can be made 
arbitrarily small in the $C^1$-topology. 
Hence we may assume w.l.o.g.\ that $U_1$ and $U_2$ 
contain small Euclidean $n$-balls of radius $R>0$. Denote the centers of 
these balls by $p_i$.  

Let $0<t_2<\min\{R, 2^{-4}\}$. Put $t_{-2}:=2^{-9}\cdot t^{16}_2$. Choose 
a smooth function $\rho:\RR\rightarrow\RR$ with
\begin{itemize}
\item[a)]
$\rho(t)=|t| \mbox{ for } |t|\ge t_{-2}$
\item[b)]
$0<\rho(t)\le t_{-2} \mbox{ for } |t|\le t_{-2}$
\item[c)]
$|\dot{\rho} (t)| \le 1 \mbox{ for all } t$.
\end{itemize}

\begin{center} 
\begin{pspicture}(-15,-3)(15,18)
\psline{->}(-10,0)(10,0)
\qdisk(-5,0){2pt}
\rput(-5,-1){$-t_{-2}$}
\qdisk(5,0){2pt}
\rput(5,-1){$t_{-2}$}
\qdisk(0,0){2pt}
\rput(0,-1){$0$}
\uput[r](10,0){$t$}
\psline{->}(0,0)(0,15)
\psline[linewidth=2pt](5,5)(10,10)
\psline[linewidth=2pt](-5,5)(-10,10)
\psline[linestyle=dashed](5,5)(0,0)
\psline[linestyle=dashed](-5,5)(0,0)
\psbezier[linewidth=2pt](-5,5)(-2.5,2)(2.5,2)(5,5)
\rput(6,7){$\rho$}

\end{pspicture}

\begin{center}
{\bf Fig. 2}
\end{center}
\end{center}

\vspace{1cm}
\noindent
We construct a metric on $M_1\# M_2$ as follows. 
On $M_i - B(p_i, t_{-2})$ take the metric of $M_i$. Recall that this 
metric is Euclidean on $B(p_i, R)$. 
Remove the two balls $B(p_1, t_{-2}) \subset M_1$ and $B(p_2, t_{-2}) 
\subset M_2$ and replace them by the cylinder
$$
[-t_{-2}, t_{-2}] \times S^{n-1}
$$
with the warped product metric
$$
ds^2 = dt^2 + \rho(t)^2 d\sigma^2
$$
where $d\sigma^2$ is the standard metric on $S^{n-1}$ of constant 
sectional curvature 1. This yields a smooth Riemannian metric $g_{t_2}$ 
on $M_1\# M_2$. 

Let $D_{t_2}$ be the corresponding operator over $X=M_1\# M_2$ of Dirac 
type over $[-t_{+2}, t_{+2}] \times S^n$ with respect to the metric 
$g_{t_2}$. 
For any operator $D$, $a \leq b$, denote by $E_{[a, b]}(D)$ the sum of all 
eigenspaces for the eigenvalues $\lambda\in[a, b]$. The Gluing Theorem 
follows readily from the following

\vspace{1cm}
CLAIM: Given $\epsilon>0$ and $\Lambda>0$ there exists 
$\delta=\delta(\epsilon, \Lambda, 
k)>0$ where $k$ is the total number of eigenvalues of $D_1$ and $D_2$ in 
the interval $(-\Lambda, \Lambda)$ such that for 
$\lambda\in(-\Lambda + 2\epsilon , \Lambda - 2\epsilon )$ and 
$0<t_2<\min\{\delta, R\}$ we have
$$
\dm E_{\{\lambda\}} (D_1) + \dm E_{\{\lambda\}} (D_2) \le 
\dm E_{[\lambda - \epsilon, \lambda + \epsilon]} (D_{t_2}) 
$$
$$
\le 
\dm E_{[\lambda - 2\epsilon, \lambda + 2\epsilon]} (D_1) + 
\dm E_{[\lambda - 2\epsilon, \lambda + 2\epsilon]} (D_2).
$$
In fact, it will turn out that 
$\delta := \min\left\{
\frac{1}{100\Lambda^2}, 2^{-4}, 
\frac{1}{2\Lambda}, \frac{1}{2(k+1)}, 2^{-17} \cdot\frac{\epsilon^2}{(k+1)^2}
\right\}$ 
does the job. 

\vspace{1cm}
Why does the Claim imply the Gluing Theorem ?
Let $\Lambda > 0$ and $\epsilon > 0$.
To see that the Gluing Theorem holds for the metric $g_{t_2}$ if 
$0<t_2<\min\{\delta, R\}$ let us assume w.l.o.g. that $\epsilon$ is so
small that any two distinct eigenvalues of $D_1 \oplus D_2$ in
$(-\Lambda , \Lambda )$ have distance at least $4\epsilon$ from one
another and that
they have distance at least $2\epsilon$ from $-\Lambda$ and from $\Lambda$.

If $\lambda \in (-\Lambda, \Lambda )$ is an eigenvalue of $D_1 \oplus
D_2$, then the Claim yields
$$\dm E_{\{\lambda\}} (D_1) + \dm E_{\{\lambda\}} (D_2) = 
\dm E_{[\lambda - \epsilon, \lambda + \epsilon]} (D_{t_2}) 
$$
Thus for any eigenvalue of $D_1 \oplus D_2$ of multiplicity $m$
there are exactly $m$ eigenvalues of
$D_{t_2}$ at distance at most $\epsilon$.

Are there further eigenvalues of $D_{t_2}$ in $(-\Lambda, \Lambda )$ ?
Since such an eigenvalue $\mu$ would have to have distance more than
$\epsilon$ from all eigenvalues of $D_1$ and $D_2$ we could find 
$\lambda \in (-\Lambda, \Lambda )$ such that $[\lambda -\epsilon,
\lambda +\epsilon ]$ contains $\mu$ but $[\lambda -2\epsilon,
\lambda +2\epsilon ]$ does not contain any eigenvalue of $D_1$ or $D_2$.
Then again by the Claim
\begin{eqnarray*}
1 & \leq & \dm E_{[\lambda - \epsilon, \lambda + \epsilon]} (D_{t_2}) \\
&\leq& \dm E_{[\lambda - 2\epsilon, \lambda + 2\epsilon]} (D_1) + 
\dm E_{[\lambda - 2\epsilon, \lambda + 2\epsilon]} (D_2) \\
&=& 0,
\end{eqnarray*}
a contradiction.
This proves the Gluing Theorem.

\vspace{1cm}
\noindent
To prove the Claim we set $t_1:=\frac{1}{2}\cdot t^4_2, 
t_{-1}:=\frac{1}{2}t_1$. 

\begin{center}
\begin{pspicture}(-6,-6)(6,6)

\psellipse[fillstyle=solid,fillcolor=lightgray](0,0)(8,4)
\psellipse[fillstyle=solid,fillcolor=darkgray](0,0)(6,3)
\psellipse[fillstyle=solid,fillcolor=lightgray](0,0)(5,2.5) 
\psellipse[fillstyle=solid,fillcolor=white](0,0)(3,1.5)

\pscustom[linecolor=white,fillstyle=solid,fillcolor=white]{
  \pscurve(2,0)(1,0.5)(0.5,2)(0.3,4.5)
  \pscurve[liftpen=1](-0.3,4.5)(-0.5,2)(-1,0.5)(-2,0)}
\pscurve(2,0)(1,0.5)(0.5,2)(0.3,4.5)
\pscurve(-2,0)(-1,0.5)(-0.5,2)(-0.3,4.5)

\rput(9,0.5){\psframebox*[framearc=0.5]{$t_2$}}
\rput(6,0.5){\psframebox*[framearc=0.5]{$t_1$}}
\rput(4,-1){\psframebox*[framearc=0.5]{$t_{-1}$}}
\rput(3,0.5){\psframebox*[framearc=0.5]{$t_{-2}$}}

\end{pspicture}
\begin{center}
{\bf Fig. 3}
\end{center}
\end{center}

\vspace{1cm}
\noindent
We choose smooth cut-off functions $\chi_i:M_i\rightarrow\RR$ such that
\begin{itemize}
\item[a)]
$0\le\chi_i\le 1$
\item[b)]
$\chi_i\equiv 1 \mbox{ on } M_i - B(p_i, t_{1})$
\item[c)]
$\chi_i\equiv 0 \mbox{ on }  B(p_i, t_{-1})$
\item[d)]
$|\nabla\chi_i| \le \frac{2}{t_1 - t_{-1}} = \frac{4}{t_1}$.
\end{itemize}

\vspace{1cm}
{\bf Proof of inequality} 
$$\dm E_{\{\lambda\}} (D_1) + \dm E_{\{\lambda\}} (D_2) \le 
\dm E_{[\lambda - \epsilon, \lambda + \epsilon]} (D_{t_2}).$$

\vspace{0.5cm}
\noindent
The idea of proof of this inequality is as follows.
Pick an eigenspinor $\sigma$ for the eigenvalue $\lambda$ on $M_1$, say.
We multiply $\sigma$ by the cut-off function $\chi_1$ to obtain a spinor
$\tilde{\sigma}$ on $X$.
Of course, $\tilde{\sigma}$ is no longer an eigenspinor.
We use $\tilde{\sigma}$ as a test function and plug it into the
Rayleigh quotient for $D_{t_2} - \lambda$ and have to show that this
Rayleigh quotient is bounded by $\epsilon$.
It then follows that $D_{t_2}$ has an eigenvalue in the 
$\epsilon$-neighborhood of $\lambda$.

Multiplying by $\chi_1$ makes the $L^2$-norm smaller, thus the denominator
of the Rayleigh quotient becomes smaller.
We have to show that we do not loose too much $L^2$-norm.
Inequality (\ref{g2a2}) says that we loose at most half of the
$L^2$-norm.
This is not a very sharp estimate but sufficient for our purposes.

More seriously, in the numerator of the Rayleigh quotient there appears
a very large error term involving the gradient of the cut-off function.
But the support of the cut-off function is contained in the annulus 
$Z^n_{t_{-1},t_1}$ (dark grey in Figure 3).
We show that the $L^2$-norm of $\sigma$ in this annulus is so small
compared to its whole $L^2$-norm that it overcompensates for the error
term, see inequality (\ref{g2a3}).

All estimates in this section use the fact that the $L^2$-norm of eigenspinors
on a Euclidean annulus is not arbitrarily distributed but tends to
cumulate near the boundary of the annulus.
In Figure 3 this means that the dark grey inner annulus carries very little
$L^2$-norm compared to the light grey region.
These ``a-priori estimates'' are proved in Section 3, see Proposition 3.2
and corollaries.

\vspace{0.5cm}
To prove the inequality we define linear maps $E_{\{\lambda\}} (D_i) 
\rightarrow L^2 (\Sigma X)$ by 
$\sigma\rightarrow\tilde{\sigma}:=\chi_i\cdot\sigma$. 
We regard $\tilde{\sigma}$ as a spinor over $X$ by extending it by 0 in
the obvious manner.
By the 
unique-continuation property for Dirac operators \cite{ar,bw} these 
maps are $1-1$. 
Since $\chi_1$ and $\chi_2$ have disjoint supports in $X$, the images 
$\tilde{E}_{\{\lambda\}} (D_1)$ and $\tilde{E}_{\{\lambda\}} (D_2)$ of 
$E_{\{\lambda\}} (D_1)$ and $E_{\{\lambda\}} (D_2)$ in 
$ L^2 (\Sigma X)$ are $L^2$-orthogonal. In particular,
$$
\dm\left(\tilde{E}_{\{\lambda\}}(D_1)+\tilde{E}_{\{\lambda\}}(D_2)\right)
= 
\dm E_{\{\lambda\}} (D_1) + \dm E_{\{\lambda\}} (D_2).
$$
We will show that the Rayleigh quotient 
of $(D_{t_2}-\lambda)^2$ is bounded by $\epsilon^2$ on 
$\tilde{E}_{\{\lambda\}} (D_1) \oplus \tilde{E}_{\{\lambda\}} (D_2)$,
provided $t_2 < \delta$. 
This obviously implies the inequality in question. 
Since $\tilde{E}_{\{\lambda\}} (D_1)$ and $\tilde{E}_{\{\lambda\}} (D_2)$ 
are orthogonal we may look at both spaces separately. 

Let $\sigma\in E_{\{\lambda\}} (D_1)$, 
$\tilde{\sigma}=\chi_1\cdot\sigma\in\tilde{E}_{\{\lambda\}} (D_1)$. 
Then
\begin{eqnarray}
\frac{\|(D_{t_2}-\lambda)\tilde{\sigma}\|^2_{L^2(X)}}
{\|\tilde{\sigma}\|^2_{L^2(X)}}
&=&
\frac{\|\nabla\chi_1\cdot\sigma+\chi_1\cdot D_1\sigma-\lambda
\cdot\chi_1\cdot\sigma\|^2_{L^2(X)}}
{\|\tilde{\sigma}\|^2_{L^2(X)}} \nonumber \\
&=&
\frac{\|\nabla\chi_1\cdot\sigma\|^2_{L^2(X)}}
{\|\tilde{\sigma}\|^2_{L^2(X)}}.
\label{g2a1}
\end{eqnarray}
By Corollary 2 to Proposition 3.2 (see next section) we know
\begin{eqnarray*}
\|\sigma\|^2_{L^2(Z^n_{0, t_1})}
&\le&
2^9 \cdot t^2_1 \cdot t_2 \cdot \|\sigma\|^2_{L^2(Z^n_{0, t_2})} \\
&\le&
2^9 \cdot 2^{-12} \cdot \|\sigma\|^2_{L^2(M_1)} 
\mbox{\hspace{3.7cm} since }t_1 \le t_2 \le \delta \le 2^{-4}\\
&\le&
\frac{1}{2} \cdot \|\sigma\|^2_{L^2(M_1)}.
\end{eqnarray*}
Hence 
\begin{equation}
\|\tilde{\sigma}\|^2_{L^2(X)}
\ge\|\sigma\|^2_{L^2(M_1)}-\|\sigma\|^2_{L^2(Z^n_{0, t_1})} 
\ge\frac{1}{2}\|\sigma\|^2_{L^2(M_1)}.
\label{g2a2}
\end{equation}

\noindent
Since $t_2 < \delta \leq \frac{1}{100\Lambda^2}$ we can apply Corollary
1 to Proposition 3.2 and obtain, using property d) of $\chi_1$
\begin{eqnarray}
\|\nabla\chi_1\cdot\sigma\|^2_{L^2(X)} 
&=&
\|\nabla\chi_1\cdot\sigma\|^2_{L^2(Z^n_{t_{-1}, t_1})} \nonumber \\
&\le&
\|\nabla\chi_1\|^2_{L^{\infty}} \cdot \|\sigma\|^2_{L^2(Z^n_{t_{-1}, t_1})} 
\nonumber \\
&\le&
\frac{16}{t^2_1} \cdot 2^7 \cdot t^2_1 \cdot t_2 \cdot 
\|\sigma\|^2_{L^2(Z^n_{t_{-2}, t_2})} \nonumber \\
&\le&
2^{11} \cdot t_2 \cdot \|\sigma\|^2_{L^2(M_1)}.
\label{g2a3}
\end{eqnarray}
By (\ref{g2a1}), (\ref{g2a2}), and (\ref{g2a3}) we can estimate the
Rayleigh quotient
$$\frac{\|(D_{t_2}-\lambda)\tilde{\sigma}\|^2_{L^2(X)}}
{\|\tilde{\sigma}\|^2_{L^2(X)}} =
\frac{\|\nabla\chi_1\cdot\sigma\|^2_{L^2(X)}}
{\|\tilde{\sigma}\|^2_{L^2(X)}}\le 2^{12} 
\cdot t_2 \le \epsilon^2$$ 
since $t_2 < \delta < 2^{-12}\cdot\epsilon^2$.

\vspace{1cm}
{\bf Proof of inequality} 
$$\dm E_{[\lambda - \epsilon, \lambda + \epsilon]} (D_{t_2}) \le 
\dm E_{[\lambda - 2\epsilon, \lambda + 2\epsilon]} (D_1) + 
\dm E_{[\lambda - 2\epsilon, \lambda + 2\epsilon]} (D_2).$$

\vspace{0.5cm}
\noindent 
This inequality is slightly more difficult to prove than the previous one.
One additional difficulty comes from the fact that this time we do not
start with a single eigenspinor but with a linear combination of 
eigenspinors for different eigenvalues.
This is where the dependence of $\delta$ on the total number $k$ of
eigenvalues of $D_1$ and $D_2$ in the interval $(-\Lambda ,\Lambda )$
comes into the game.

Secondly, since we start with eigenspinors on $X$ we also have to control
their $L^2$-norm on the connecting neck between $M_1$ and $M_2$.
Hence we need estimates for eigenspinors on cylindrical manifolds with
certain warped product metrics.
This is provided by Propsition 3.3.
The rest of the proof is analogous to that of the previous inequality.

\vspace{0.5cm}
For $-t_2\le\mu<\nu\le t_2$ denote the connecting cylinder 
$[\mu, \nu] \times S^{n-1}\subset X$ with metric $ds^2 = dt^2 + \rho^2 
d\sigma^2$ by $Z_{\mu, \nu}$. 

We assume $\Lambda\cdot t_2\le\Lambda\cdot\delta\le\frac{1}{2}$ so that for 
$|\lambda|\le\Lambda$: 
$$
|\lambda| \cdot |\rho| + \frac{1}{2} |\dot{\rho}| \le 
\Lambda\cdot t_2 + \frac{1}{2} \cdot 1 \le 1.
$$
Hence we can apply Proposition 3.3 and obtain for any eigenspinor 
$\sigma$ on $Z_{-t_2, t_2}$ for the eigenvalue $\lambda$
\begin{eqnarray*}
\|\sigma\|^2_{L^2(Z_{-t_1, t_1})} 
&\le&
\frac{t_1 - (-t_1)}{2(t_2 - t_1)} \cdot 
\|\sigma\|^2_{L^2(Z_{-t_2, t_2})} \\
&=&
\frac{t^4_2}{2(1-\frac{1}{2}t^3_2)t_2} \cdot 
\|\sigma\|^2_{L^2(Z_{-t_2, t_2})} \\
&\le&
t^3_2 \cdot \|\sigma\|^2_{L^2(Z_{-t_2, t_2})}.
\end{eqnarray*}
If $\sigma=\Sigma\sigma_j$ is a sum of at most $k+1$ 
$D_{t_2}$-eigenspinors on $X$,
pairwise $L^2(X)$-orthogonal, then we conclude
\begin{eqnarray*}
\|\sigma\|_{L^2(Z_{-t_1, t_1})}
&\le&
\sum\limits_j \|\sigma_j\|_{L^2(Z_{-t_1, t_1})} \\
&\le&
t^{3/2}_2 \cdot\sum\limits_j \|\sigma_j\|_{L^2(Z_{-t_2, t_2})} \\
&\le&
t^{3/2}_2 \cdot\sum\limits_j \|\sigma_j\|_{L^2(X)} \\
&\le&
t_2 \cdot (k+1) \cdot \|\sigma\|_{L^2(X)}.
\end{eqnarray*}
Since we assume $t_2 < \delta\le\frac{1}{2(k+1)}$, we conclude
\begin{equation}
\|\sigma\|_{L^2(Z_{-t_1, t_1})} \le \frac{1}{2} 
\|\sigma\|_{L^2(X)}.
\label{g2.1}
\end{equation}
Put $l:=\dm E_{[\lambda - 2\epsilon, \lambda + 2\epsilon]} (D_1) + 
\dm E_{[\lambda - 2\epsilon, \lambda + 2\epsilon]} (D_2)$. 
We will derive a contradiction from the assumption
$$
\dm E_{[\lambda - \epsilon, \lambda + \epsilon]} (D_{t_2}) \ge l+1.
$$
Let $V$ be spanned by $l+1$ linearly independent $D_{t_2}$-eigenspinors 
for eigenvalues in $[\lambda - \epsilon, \lambda + \epsilon]$. 
Note that $l+1\le k+1$. The linear map
\begin{eqnarray*}
V &\rightarrow& L^2 (\Sigma M_1) \oplus  L^2 (\Sigma M_2) \\
\sigma &\rightarrow& (\chi_1\cdot\sigma, \chi_2\cdot\sigma)
\end{eqnarray*}
is again $1-1$ by the unique-continuation property. 
We will show that the Rayleigh quotient of 
$(D_1 - \lambda)^2 \oplus (D_2 - \lambda)^2$ is bounded by 
$(2\epsilon )^2$ on the image of $V$. 
This would imply that $D_1\oplus D_2$ has at least $l+1$ eigenvalues in 
$[\lambda - 2\epsilon, \lambda + 2\epsilon]$, a contradiction.
\begin{eqnarray}
\|(\chi_1\cdot\sigma, \chi_2\cdot\sigma)\|^2_{L^2(M_1) \oplus L^2(M_2)} 
&=&
\|\chi_1\cdot\sigma\|^2_{L^2(M_1)} + 
\|\chi_2\cdot\sigma\|^2_{L^2(M_2)} \nonumber \\
&\ge&
\|\sigma\|^2_{L^2(X)} - \|\sigma\|^2_{L^2(Z_{-t_1, t_1})} \nonumber \\
&\ge&
\frac{1}{2} \|\sigma\|^2_{L^2(X)}
\label{g2b1}
\end{eqnarray}
by (\ref{g2.1}).
This implies for the Rayleigh quotient
$$
\frac{
\|(D_1-\lambda)\oplus(D_2-\lambda) (\chi_1\sigma,
\chi_2\sigma)\|^2_{L^2(M_1)\oplus L^2(M_2)}
}
{
\|(\chi_1\sigma, \chi_2\sigma)\|^2_{L^2(M_1)\oplus L^2(M_2)}
} =
$$
\begin{eqnarray}
&=& 
\frac{
\|(D_1-\lambda)(\chi_1\sigma)\|^2_{L^2(M_1)} +
\|(D_2-\lambda)(\chi_2\sigma)\|^2_{L^2(M_2)}
}
{
\|(\chi_1\sigma, \chi_2\sigma)\|^2_{L^2(M_1)\oplus L^2(M_2)}
} \nonumber \\
&\leq&
2\cdot\frac{
\|(D_1-\lambda)(\chi_1\sigma)\|^2_{L^2(M_1)} +
\|(D_2-\lambda)(\chi_2\sigma)\|^2_{L^2(M_2)}
}
{
\|\sigma\|^2_{L^2(X)}
}
\label{2xxx}
\end{eqnarray}
By property d) of $\chi_i$ and Corollary 1 to Proposition 3.2 we know
\begin{eqnarray}
\|\nabla\chi_i\cdot\sigma\|_{L^2(M_i)}
&=&
\|\nabla\chi_i\cdot\sigma\|_{L^2(Z^n_{t_{-1}, t_1})} \nonumber \\
&\le& 
\frac{4}{t_1} \cdot\|\sigma\|_{L^2(Z^n_{t_{-1}, t_1})} \nonumber \\
&\le& 
\frac{4}{t_1} \cdot\sqrt{2^7\cdot t^2_1 \cdot t_2} \cdot 
\sum\limits_j 
\|\sigma_j\|_{L^2(Z^n_{t_{-2}, t_2})} \nonumber \\
&\le&
\sqrt{2^{11}\cdot t_2} \cdot (k+1) \cdot \|\sigma\|_{L^2 (X)}
\nonumber \\
&\le& 
\frac{\epsilon}{8}\cdot\|\sigma\|_{L^2(X)}
\label{g2b2}
\end{eqnarray}
because $t_2 < \delta \leq 2^{-17}\cdot\frac{\epsilon^2}{(k+1)^2}$.
Hence, by (\ref{g2b2})
$$
\frac{
\|(D_i-\lambda)(\chi_i\sigma)\|_{L^2(M_i)}
}
{
\|\sigma\|_{L^2(X)}
}
=
\frac{\|\chi_i (D_i-\lambda)\sigma + \nabla\chi_i\sigma\|_{L^2(M_i)}}
{\|\sigma\|_{L^2(X)}}
$$
\begin{eqnarray}
&\leq&
\frac{\|\chi_i (D_i-\lambda)\sigma\|_{L^2(M_i)} + 
\|\nabla\chi_i\sigma\|_{L^2(M_i)}}
{\|\sigma\|_{L^2(X)}} \nonumber \\
&\leq&
\frac{
\|(D_{t_2}-\lambda)(\sigma)\|_{L^2(M_i-B(p_i,t_{-1}))}}
{\|\sigma\|_{L^2(X)}}
+ \frac{\epsilon}{8}
\label{2xxy}
\end{eqnarray}
Combining (\ref{2xxx}) and (\ref{2xxy}) we obtain the desired estimate
for the Rayleight quotient:
$$
\frac{
\|(D_1-\lambda)\oplus(D_2-\lambda) (\chi_1\sigma,
\chi_2\sigma)\|^2_{L^2(M_1)\oplus L^2(M_2)}
}
{
\|(\chi_1\sigma, \chi_2\sigma)\|^2_{L^2(M_1)\oplus L^2(M_2)}
} \leq
$$
\begin{eqnarray*}
&\leq& 2 \cdot \left\{ \left[\frac{
\|(D_{t_2}-\lambda)(\sigma)\|_{L^2(M_1-B(p_1,t_{-1}))}}
{\|\sigma\|_{L^2(X)}}
+ \frac{\epsilon}{8} \right]^2 + \right. \\
&&
+ \left.\left[\frac{
\|(D_{t_2}-\lambda)(\sigma)\|_{L^2(M_2-B(p_2,t_{-1}))}}
{\|\sigma\|_{L^2(X)}}
+ \frac{\epsilon}{8} \right]^2
\right\} \\
&=&
2\cdot\left\{
\frac{
\|(D_{t_2}-\lambda)(\sigma)\|^2_{L^2(X-Z_{-t_{-1},t_{-1}})}}
{\|\sigma\|^2_{L^2(X)}} + 2\cdot \left(\frac{\epsilon}{8}\right)^2 
+ \right. \\
&&
+ \left. 2\cdot\frac{\epsilon}{8}\cdot\frac{
\|(D_{t_2}-\lambda)(\sigma)\|_{L^2(M_1-B(p_1,t_{-1}))} + 
\|(D_{t_2}-\lambda)(\sigma)\|_{L^2(M_2-B(p_2,t_{-1}))}}
{\|\sigma\|_{L^2(X)}}
\right\} \\
&\leq&
2\cdot\{\epsilon^2 + 2\cdot(\frac{\epsilon}{8})^2 + 
2\cdot\frac{\epsilon}{8}\cdot2\epsilon \} \\
&<& (2\epsilon)^2.
\end{eqnarray*} 

\noindent
This finishes the proof of the Claim and of the Gluing Theorem. \hfill$\Box$

\vspace{2cm}
\section{The Estimates}
In this section we will derive the a-priori estimates on the distribution 
of the $L^2$-norm of eigenspinors on Euclidean annuli which have been 
used in the previous section to prove the Gluing Theorem. 
Roughly speaking, what we show is the fact that the $L^2$-norm of
an eigenspinor on a Euclidean annulus cumulates in the region near 
the boundary.
In Figure 3 this means that there is only very little $L^2$-norm
on the dark grey inner annulus $Z^n_{t_{-1},t_1}$ compared to the
$L^2$-norm on the whole annulus $Z^n_{t_{-2},t_2}$, compare Corollary
1 to Propostion 3.2.
To prove this we regard a Euclidean annulus as a warped product of an
interval with a sphere and decompose the eigenspinor with respect
to an eigenbasis on the sphere.
Then the eigenspinor equation translates into ordinary differential 
equations on the coefficients.

To controle the coefficients 
we prove an estimate which allows us to compare the unknown solution of 
a linear ordinary differential equation with certain ``almost solutions''. 
In the proof of Proposition 3.2 we will deal with linear
ordinary differential equations which we cannot explicitly solve.
But we will be able to guess the right ``almost solutions'' and thus
obtain information on the unknown solutions.
Proposition 3.1 and its proof are similar to the well-known Gronwall
lemma.

\vspace{1cm}
{\bf Proposition 3.1.} {\it Let $I\subset\RR$ be an interval}, $t_0 \in 
I$, {\it let} $A:I\rightarrow Mat (n\times n, \CC)$ {\it be a continuous 
mapping into the complex $n\times n$-matrices}. {\it Let $u$ be a solution of}
$$
\dot{u} (t)=A(t)u(t).
$$
{\it Moreover, let $v:I\rightarrow\CC^n$ be a continuously differentiable 
function such that}
\begin{itemize}
\item[(i)]
$v(t_0)=u(t_0)$
\item[(ii)]
$|\dot{v}(t) - A(t)v(t)|\le\delta(t)$
\end{itemize}
{\it where} $\delta:I\rightarrow[0, \infty)$ {\it is a continuous 
function}. 

{\it Then the following estimate holds}
$$
|u(t)-v(t)| \le 
\left| \; 
\int\limits^t_{t_0} \delta(s)\cdot e^{\| A\|_{\infty}\cdot|t-s|} ds
\right|
$$
{\it where} 
$\| A\|_{\infty} = 
\sup\limits_{t\in I} |A(t)|_{\mbox{\scriptsize Op}}$ {\it and} 
$|A(t)|_{\mbox{\scriptsize Op}}$ {\it is the operator norm}.

\vspace{1cm}
{\bf Proof.} It is sufficient to prove the claim for $t\ge t_0$. 
\newline
Put $w(t):=u(t)-v(t)$. Then
$$
\dot{w}(t)=A(t)\cdot w(t)+\rho(t)
$$
where $\rho(t)=A(t)\cdot v(t)-\dot{v}(t)$. Thus
\begin{eqnarray}
|w(t)|&=&\left|
\int\limits^t_{t_0} (A(s)\cdot w(s) + \rho(s)) ds
\right| \nonumber \\  
&\le&
\| A\|_{\infty} \cdot \int\limits^t_{t_0} |w(s)| ds +
\int\limits^t_{t_0} \delta(s) ds. \label{g3.3}
\end{eqnarray}
For $\varphi(t):=\| A\|_{\infty} \cdot \int\limits_{t_0}^t 
|w(s)| ds + \int\limits_{t_0}^t \delta(s) ds$ 
\newline
we have
$$
\varphi(t_0) = 0 \\
$$
and
\begin{equation}
\dot{\varphi} (t) = \| A\|_{\infty} 
|w(t)| + \delta(t).
\label{g3.4}
\end{equation}
Hence by (\ref{g3.3}) and (\ref{g3.4})
\begin{eqnarray*}
0&\le& \| A\|_{\infty} \cdot 
\left(\varphi(t)-|w(t)|\right) \\
&=& \| A\|_{\infty} \cdot \varphi(t) - \dot{\varphi}(t) + \delta(t).
\end{eqnarray*}
Multiplying by $e^{- \| A\|_{\infty} \cdot (t - t_0)}$ we 
obtain 
$$
0\le - \frac{d}{dt}\left(
e^{- \| A\|_{\infty} \cdot (t - t_0)} \cdot \varphi(t)
\right)
+ e^{-\| A\|_{\infty} \cdot (t - t_0)} \cdot \delta(t).
$$
Integration yields
$$
e^{- \| A\|_{\infty} \cdot (t - t_0)} \cdot \varphi(t) 
\le 
\int\limits^t_{t_0} \left(e^{- \| A\|_{\infty} 
\cdot (s - t_0)} \cdot \delta(s)\right) ds
$$
which implies again using (\ref{g3.3})
\begin{eqnarray*}
|w(t)| &\le& \varphi(t) \\
 &\le& \int\limits^t_{t_0} e^{\| A\|_{\infty} (t-s)} \cdot 
 \delta(s) ds. \hspace{6.7cm}\Box
\end{eqnarray*}

\vspace{1cm}
\noindent
For $0<a<b$ let $Z^n_{a, b} := \{ x\in\RR^n \mid a\le|x|\le b\}$ be the 
corresponding {\it Euclidean annulus}. As a Riemannian manifold we can also 
describe $Z^n_{a, b}$ as a warped product
$$
Z^n_{a, b} = [a, b] \times S^{n-1}
$$
with the metric
$$
ds^2 (t, y)= dt^2 + t^2 d\sigma^2 (y)
$$
where $d\sigma^2$ is the standard metric on $S^{n-1}$ of constant 
sectional curvature 1.

Restriction of the spinor bundle $\Sigma\RR^n$  of $\RR^n$ to $S^{n-1}$ 
yields the spinor bundle of $S^{n-1}$ if $n$ is odd. If $n$ is even one 
obtains the direct sum of two copies of the spinor bundle on $S^{n-1}$. 
If $n$ is odd let $\tilde{D}$ be the Dirac operator of $S^{n-1}$,
if $n$ is even let $\tilde{D}$ be the sum of the Dirac operator of 
$S^{n-1}$ and of its negative.

We choose an orthonormal basis of $L^2 (\Sigma\RR^n|_{S^{n-1}})$ 
consisting of eigenspinors $\sigma_j$ of $\tilde{D}$. Clifford  
multiplication by the normal vector field $\frac{\partial}{\partial t}$ 
anticommutes with $\tilde{D}$. Thus if $\sigma_j$ is an eigenspinor for 
the eigenvalue $\mu_j$, then $\frac{\partial}{\partial t} \cdot \sigma_j$ is 
an eigenspinor for $-\mu_j$. Hence we may assume 
$\sigma_{-j} = \frac{\partial}{\partial t} \cdot \sigma_j$, $j\in\NN$, 
$0<\mu_1\le\mu_2\le\mu_2\le\ldots\rightarrow +\infty$. 
Parallel translation along the $t$-lines yields smooth sections $(t, y) 
\mapsto \sigma_j(t, y)$ on $(0, \infty) \times S^{n-1} = \RR^n - \{0\}$, 
$j\in\ZZ^* = \ZZ - \{0\} = \NN \cup - \NN$. If $\sigma$ is a smooth 
spinor field on $Z^n_{a, b}$, then we can write
$$
\sigma(t, y)=\sum\limits_{j\in\ZZ^*} \beta_j(t) \cdot t^{-\frac{n-1}{2}} 
\cdot \sigma_j (t, y)
$$
where $\beta_j$ are smooth functions, $\beta_j : [a, b] \rightarrow\CC$.
The normalization factor $t^{-\frac{n-1}{2}}$ is chosen such that the
subsequent formulas become very simple. 
If $D$ is the Dirac operator on $Z^n_{a, b}$, then the eigenspinor equation
$$
D\sigma=\lambda\cdot\sigma
$$
is equivalent to 
\begin{equation}
\dot{B}_j (t) = 
\left(
\matrix{\frac{\mu_j}{t} & -\lambda\cr
\cr
\lambda & -\frac{\mu_j}{t}}
\right) 
\cdot B_j(t)
\label{g3.1}
\end{equation}
where $B_j(t) = \left(\matrix{\beta_{-j}(t) \cr \beta_j (t)} 
\right)$ for all $j \in \NN$. 
The $L^2$-norm of $\sigma$ is given by 
$$
\|\sigma\|^2_{L^2 (Z^n_{a, b})} = \sum\limits^{\infty}_{j=1} 
\int\limits^b_a |B_j(t)|^2 dt.
$$
See \cite{b4} for the details. 
Now we come to the main estimate of this section.

\vspace{1cm}
{\bf Proposition 3.2.} {\it Let} $0<t_{-2}<t_{-1}<t_1 < t_2\le 1$, 
{\it let} $\mu\ge 1$, {\it let} $\lambda\in\RR$. 
{\it If} $|\lambda| t^{1/2}_2 \le \frac{1}{10}$, $t_2\ge 2t_1$, 
$t_{-1}\ge 2t_{-2}$, {\it and} $t^6_1\le t_{-2}$, 
{\it then every solution $B$ of}
$$
\dot{B}(t)=\left(
\matrix{
\frac{\mu}{t} & - \lambda\cr
\cr
\lambda&-\frac{\mu}{t}
}
\right) \cdot B(t) \qquad \label{(a)}
$$
{\it satisfies} 
$$
\frac{\int\limits^{t_1}_{t_{-1}}|B(t)|^2 dt}
{\int\limits^{t_2}_{t_{-2}}|B(t)|^2 dt}
\le 
2^6 \cdot \max\left\{
3\cdot \left(\frac{t_1}{t_2}\right)^{2\mu+1}, 
\frac{t_1}{t_{-1}}\cdot\left(\frac{t_{-2}}{t_{-1}}\right)^{2\mu-1} 
\right\}
$$

\vspace{1cm}
{\bf Corollary 1.} {\it Let} $n\geq 3$, {\it let} $0<t_2<2^{-4}$, 
{\it let} $\Lambda>0$. 
{\it Put} $t_1:=\frac{1}{2}\cdot t^4_2$, $t_{-1}:=\frac{1}{2}t_1$, 
$t_{-2}:=\frac{1}{2}t^4_{-1}$ {\it Then every eigenspinor} $\sigma$ 
{\it for the eigenvalue} $\lambda$ {\it on} $Z^n_{t_{-2}, t_2}$ {\it satisfies}
$$
\frac{
\|\sigma\|^2_{L^2 (Z^n_{t_{-1}, t_1})}
}{
\|\sigma\|^2_{L^2 (Z^n_{t_{-2}, t_2})}
}
\le2^7\cdot t^2_1\cdot t_2,
$$
{\it provided} $|\lambda| \le \Lambda$ {\it and} 
$\Lambda \cdot t^{1/2}_2 \le \frac{1}{10}$.

\vspace{1cm}
{\bf Proof of Corollary 1.} With our definition of $t_1$, $t_{-1}$, and 
$t_{-2}$ the assumptions on $t_j$ in Proposition 3.2 are trivially 
satisfied. 
Moreover, it is well-known that the absolute value of all 
Dirac eigenvalues of $S^{n-1}$ is at least $1$ if $n\ge 3$, compare e.g. 
\cite{b1}. Hence
\begin{eqnarray*}
\|\sigma\|^2_{L^2 (Z^n_{t_{-1}, t_1})} &=&
\sum\limits^{\infty}_{j=1} \; 
\int\limits^{t_1}_{t_{-1}} |B_j(t)|^2 dt\\
&\le&
2^6\cdot\sum\limits^{\infty}_{j=1}\max \left\{
3\cdot
\left( \frac{1}{2}t^3_2\right)^{2\mu_j +1}, \, 
2\cdot\left( \frac{1}{2}t^3_{-1}\right)^{2\mu_j -1}
\right\} 
\cdot \int\limits^{t_2}_{t_{-2}} |B_j(t)|^2 dt \\
&\le&
2^6 \cdot \max \left\{3 \cdot 2^{-3}\cdot t^9_2, \, t^3_{-1}
\right\} \cdot 
\|\sigma\|^2_{L^2 (Z^n_{t_{-2}, t_2})} \\
&\le&
2^7 \cdot t^2_1 \cdot t_2 \cdot 
\|\sigma\|^2_{L^2 (Z^n_{t_{-2}, t_2})}. \hspace{6.3cm}\Box
\end{eqnarray*}

\vspace{1cm}
{\bf Corollary 2.} {\it Let} $n\geq 3$, {\it let} $0<t_2<2^{-4}$, {\it let} 
$t_1 := \frac{1}{2}t^4_2$, {\it let} $\Lambda>0$. 
{\it Then every eigenspinor} $\sigma$ {\it for the eigenvalue} $\lambda$ 
{\it on the ball} $Z^n_{0, t_2}$ {\it satisfies}
$$
\frac{
\|\sigma\|^2_{L^2 (Z^n_{0, t_1})}
}{
\|\sigma\|^2_{L^2 (Z^n_{0, t_2})}
}
\le 2^9 \cdot t^2_1 \cdot t_2
$$
{\it if} $|\lambda|\le\Lambda$ {\it and} 
$\Lambda t^{1/2}_2 \le \frac{1}{10}$.

\vspace{1cm}
{\bf Proof of Corollary 2.} We apply Corollary 1 to $(k\in\NN_0)$
\begin{eqnarray*}
t_{2, k}   &:=&  2^{-k/4} \cdot t_2,    \\ 
t_{1, k}   &:=& 2^{-k}   \cdot t_1,  \\ 
t_{-1, k}  &:=&  \frac{1}{2}t_{1,k} = 2^{-k-1}   \cdot t_{1}, \\ 
t_{-2, k}  &:=&  \frac{1}{2}t^4_{-1,k} :    
\end{eqnarray*}
\begin{eqnarray*}
\|\sigma\|^2_{L^2 (Z^n_{t_{-1, k}, t_{1, k}})} 
&\le&
2^7 \cdot t^2_{1, k} \cdot t_{2, k} \cdot 
\|\sigma\|^2_{L^2 (Z^n_{t_{-2, k}, t_{2, k}})} \\
&\le&
2^7 \cdot 2^{-2k} \cdot t^2_1 \cdot t_2 \cdot 
\|\sigma\|^2_{L^2 (Z^n_{0, t_2})}.
\end{eqnarray*}
Summation over $k$ yields
$$
\|\sigma\|^2_{L^2 (Z^n_{0, t_1})} \le 2^7 \cdot \frac{4}{3} \cdot 
t^2_1 \cdot t_2 \cdot 
\|\sigma\|^2_{L^2 (Z^n_{0, t_2})} \hspace{6.1cm}\Box
$$

\vspace{1cm}
{\bf Proof of Proposition 3.2.} The substitution $t=e^{\tau}$, 
$\tilde{B}(\tau)=B(e^{\tau})$, translates our differential equation 
(\ref{g3.1}) into 
\begin{equation}
\tilde{B}'(\tau)=\left(
\matrix{\mu&-\lambda e^{\tau}\cr
\cr
\lambda e^{\tau}&-\mu}
\right)
\cdot\tilde{B}(\tau). 
\label{g3.2}
\end{equation}
Set $t_0 := \sqrt{t_1 \cdot t_{-1}}$, $t_k=e^{\tau_k}$, 
$k=-2, -1, 0, 1, 2$. Write 
$$
B(t_0)=\tilde{B}(\tau_0)={w_1\choose w_2}.
$$
We will compare the solution $\tilde{B}$ of (\ref{g3.2}) to the 
``almost solution''
\begin{eqnarray*}
\tilde{v}(\tau) 
&:=& 
\left(\matrix{w_1 \cdot e^{\mu (\tau - \tau_0)} \cr
w_2 \cdot e^{-\mu(\tau - \tau_0)}}\right), \\
\\
v(t) 
&:=& 
\left(\matrix{w_1 \cdot \left(\frac{t}{t_0}\right)^{\mu} \cr
w_2 \cdot \left(\frac{t}{t_0}\right)^{- \mu}}\right).
\end{eqnarray*}
Define
$$
A(\tau ):=\left(
\matrix{
\mu&-\lambda e^{\tau} \cr
\lambda e^{\tau}&-\mu
}
\right) .
$$
Then on $[\tau_{-2}, \tau_2]$ we have $\|A\|_{\infty} \le \mu + 
\frac{1}{10}$ because $|\lambda e^{\tau}|\le|\lambda t_2|\le
|\lambda t_2^{1/2}|\le\frac{1}{10}$. 
We compute
\begin{eqnarray*}
|\tilde{v}' - A(\tau) \cdot \tilde{v}(\tau)|^2 
&=& 
\left|
\lambda \cdot e^{\tau} \cdot \left(
\matrix{w_2 \cdot e^{-\mu (\tau - \tau_0)} \cr
-w_1 \cdot e^{\mu (\tau - \tau_0)}}
\right)
\right|^2 \\
\\
&=&
\lambda^2 \cdot e^{2\tau} \cdot 
\left(
|w_1|^2 e^{2\mu(\tau - \tau_0)} + |w_2|^2 e^{-2\mu(\tau - \tau_0)}
\right).
\end{eqnarray*}
Now we estimate the difference between the solution $\tilde{B}$ and the 
``almost solution'' $\tilde{v}$.

%\vspace{1cm}
%{\bf Case 1.} $\tau\ge\tau_0$.

%\vspace{0.5cm}
%\noindent
%Put $\delta(\tau):=\lambda\cdot e^{\tau}\cdot|B(t_0)|\cdot 
%e^{\mu(\tau-\tau_0)}$. Then
%$$
%|\tilde{v}'(\tau)-A(\tau)\cdot \tilde{v}(\tau)|\le\delta(\tau)
%$$ 
%and hence by Proposition 3.1
%\begin{eqnarray*}
%|\tilde{B}(\tau) - \tilde{v}(\tau)| 
%&\le& 
%\int\limits^{\tau}_{\tau_0} 
%\delta(\sigma) \cdot e^{\|A\|_{\infty}\cdot(\tau - \sigma)} d\sigma \\
%&\le& 
%\lambda \cdot |B(t_0)| \cdot \int\limits^{\tau}_{\tau_0} e^{\sigma} \cdot
%e^{\mu(\sigma - \tau_0)} \cdot e^{(\mu+\frac{1}{10})(\tau - \sigma)} d\sigma \\
%&=&
%\lambda \cdot e^{\mu(\tau - \tau_0)} \cdot |B(t_0)| \cdot \frac{10}{9} 
%\cdot e^{\frac{\tau}{10}} \cdot
%\left(
%e^{\frac{9}{10}\tau} - e^{\frac{9}{10}\tau_0}\right) \\
%&\le&
%\frac{10}{9} \cdot \lambda \cdot e^{\tau} \cdot |B(t_0)| \cdot 
%e^{\mu(\tau - \tau_0)} \\
%&\le&
%|B(t_0)| \cdot e^{\mu(\tau - \tau_0)}.
%\end{eqnarray*}
%Hence 
%$$
%|B(t)-v(t)| \le |B(t_0)| \cdot \left(\frac{t}{t_0}\right)^\mu .
%$$

\vspace{1cm}
%{\bf Case 2.} 
We start with the case that $\tau\in [\tau_{-2},\tau_0]$.

\vspace{0.5cm}
\noindent
We set $\delta(\tau):=|\lambda |\cdot e^{\tau}\cdot |B(t_0)|\cdot 
e^{-\mu(\tau - \tau_0)}$. 
We have
$$
|\tilde{v}'(\tau) - A(\tau) \cdot \tilde{v}(\tau)| \le \delta(\tau).
$$
Proposition 3.1. yields
\begin{eqnarray*}
|\tilde{B}(\tau) - \tilde{v}(\tau)| 
&\le& 
\int\limits^{\tau_0}_{\tau}  
\delta(\sigma) \cdot e^{\|A\|_{\infty}\cdot(\sigma - \tau)} d\sigma \\
&\le& 
|\lambda |\cdot |B(t_0)| \cdot \int\limits^{\tau_0}_{\tau} e^{\sigma} \cdot 
e^{-\mu(\sigma - \tau_0)} \cdot e^{(\mu + \frac{1}{10})(\sigma-\tau)} 
d\sigma \\ 
&=& 
|\lambda |\cdot |B(t_0)| \cdot e^{-\mu(\tau - \tau_0)} \cdot 
e^{-\frac{\tau}{10}} \cdot \int\limits^{\tau_0}_{\tau} 
e^{\frac{11}{10}\sigma} d\sigma \\ 
&\le& 
|\lambda |\cdot |B(t_0)| \cdot e^{-\mu(\tau - \tau_0)} \cdot e^{\tau_0} \cdot 
e^{\frac{1}{10}(\tau_0 - \tau)}.
\end{eqnarray*}
Thus
\begin{eqnarray}
|B(t) - v(t)| 
&\le& 
|\lambda |\cdot t_0 \cdot |B(t_0)| \cdot 
\left(\frac{t_0}{t_{-2}}\right)^{\frac{1}{10}} \cdot 
\left(\frac{t}{t_0}\right)^{-\mu} \nonumber \\
&=&
|\lambda |\cdot t^{1/2}_0 \cdot |B(t_0)| \cdot 
\left(\frac{t^6_0}{t_{-2}}\right)^{\frac{1}{10}} \cdot 
\left(\frac{t}{t_0}\right)^{-\mu} \label{3di21.3} \\
&\le&
|B(t_0)| \cdot \left(\frac{t}{t_0}\right)^{-\mu}
\label{3di21.4}
\end{eqnarray}
because $t^6_0 \le t^6_1 \le t_{-2}$ and $|\lambda |\cdot t_0^{1/2}
\leq |\lambda |\cdot t_2^{1/2} < 1$.

\vspace{1cm}
In case $\tau\in [\tau_0,\tau_2]$ similar reasonning yields
\begin{eqnarray}
|B(t)-v(t)| &\le& \frac{10}{9}\cdot |\lambda |\cdot t\cdot
|B(t_0)| \cdot \left(\frac{t}{t_0}\right)^\mu 
\label{3di21.1} \\
&\leq&
|B(t_0)| \cdot \left(\frac{t}{t_0}\right)^\mu 
\label{3di21.2}.
\end{eqnarray}

\vspace{1cm}
Now that we have some control on how much $v$ can differ from the 
solution $B$ we can look at the growth of $B$. An upper bound is now 
easily obtained.

If $t\in [ t_0,t_2 ]$ we get by (\ref{3di21.2}):
$$
|B(t)| \le |v(t)| + |B(t) - v(t)| \le 
2\cdot |B(t_0)| \cdot \left(\frac{t}{t_0}\right)^{\mu}
$$
and similarly for $t\in [t_{-2},t_0]$:
$$
|B(t)| \le 2 \cdot |B(t_0)| \cdot \left(\frac{t}{t_0}\right)^{-\mu}.
$$
>From this we get

\begin{eqnarray}
\int\limits^{t_1}_{t_{-1}} |B(t)|^2 dt 
&=&
\int\limits^{t_0}_{t_{-1}} |B(t)|^2 dt + 
\int\limits^{t_1}_{t_0}    |B(t)|^2 dt \nonumber \\
&\le&
4\cdot |B(t_0)|^2 \cdot 
\left\{
\frac{1}{2\mu-1} \cdot
\left[
\left(\frac{t_{-1}}{t_0}\right)^{-2\mu} \cdot t_{-1} -t_0
\right] \right.\nonumber \\
 & &
\hspace{2.5cm}
\left.
+\frac{1}{2\mu+1} 
\left[
\left(\frac{t_1}{t_0}\right)^{2\mu} \cdot t_1 -t_0
\right]
\right\} \nonumber \\
&\le&
\frac{8}{2\mu-1} \cdot |B(t_0)|^2 \cdot 
\left(\frac{t_1}{t_{-1}}\right)^{\mu} \cdot t_1.
\label{g3.5}
\end{eqnarray}
by the definition of $t_0$. 

To get a lower bound on $|B(t)|$ we have to distinguish two cases.

\vspace{1cm}
{\bf Case 1.} $|w_1| \ge |w_2|$.

\vspace{0.5cm}
\noindent
We know $|v(t)|^2 \ge |w_1|^2 \cdot 
\left(\frac{t}{t_0}\right)^{2\mu} \ge \frac{1}{2} |B(t_0)|^2 \cdot 
\left(\frac{t}{t_0}\right)^{2\mu}$. 
For $t\in [t_0,t_2]$ we get
\begin{eqnarray*}
\mbox{\hspace{2cm}} |B(t)| 
&\ge&
|v(t)| - |B(t) - v(t)| \\
&\ge&
|B(t_0)| \cdot \left(\frac{t}{t_0}\right)^{\mu} \cdot 
\left\{
\frac{1}{\sqrt{2}} - \frac{10}{9} \cdot |\lambda |\cdot t \right\} 
\mbox{\hspace{2.55cm}by (\ref{3di21.1})}\\
&\ge&
\frac{1}{2} \cdot |B(t_0)| \cdot \left(\frac{t}{t_0}\right)^{\mu}.
\end{eqnarray*}
Hence by (\ref{3di21.4})
\begin{eqnarray}
\int\limits^{t_2}_{t_{-2}} |B(t)|^2 dt 
&\ge&
\int\limits^{t_2}_{t_0} |B(t)|^2 dt \nonumber \\
&\ge&
\frac{1}{4} \cdot |B(t_0)|^2 \cdot \frac{1}{2\mu+1} \cdot
\left\{
\left(\frac{t_2}{t_0}\right)^{2\mu} \cdot t_2 - t_0\right\} \nonumber \\
&\ge&
\frac{1}{8} \cdot |B(t_0)|^2 \cdot \frac{1}{2\mu+1} \cdot 
\left(\frac{t_2}{t_0}\right)^{2\mu} \cdot t_2
\label{g3.6}
\end{eqnarray}
because $t_2\ge 2\cdot t_0$.
>From (\ref{g3.5}), (\ref{g3.6}), and $t_0^2 = t_1 \cdot t_{-1}$ we deduce
\begin{eqnarray*}
\frac{\int\limits^{t_1}_{t_{-1}} |B(t)|^2 dt}
{\int\limits^{t_2}_{t_{-2}} |B(t)|^2 dt} 
&\le&
\frac{\frac{8}{2\mu-1} \cdot |B(t_0)|^2 \cdot  
    \left(\frac{t_1}{t_{-1}}\right)^{\mu} \cdot t_1} 
{\frac{1}{8} \cdot |B(t_0)|^2 \cdot \frac{1}{2\mu+1} \cdot 
\left(\frac{t_2}{t_0}\right)^{2\mu} \cdot t_2}
\cdot \frac{t_1^\mu \cdot t_{-1}^\mu}{t_0^{2\mu}} \\
&=  &
2^6 \cdot \frac{2\mu+1}{2\mu-1} \cdot 
\left(\frac{t_1}{t_2}\right)^{2\mu+1} \\
&\le&
2^6 \cdot 3 \cdot \left(\frac{t_1}{t_2}\right)^{2\mu+1}.
\end{eqnarray*}

\vspace{1cm}
{\bf Case 2.} $|w_1| \le |w_2|$.

\vspace{0.5cm}
\noindent
This time we have $|v(t)|^2 \ge |w_2|^2 \cdot 
\left(\frac{t}{t_0}\right)^{-2\mu} \ge \frac{1}{2} 
|B(t_0)|^2 \cdot 
\left(\frac{t}{t_0}\right)^{-2\mu}$. 
\newline
For $t\in [t_{-2},t_0]$ we get
\begin{eqnarray*}
\mbox{\hspace{2cm}} |B(t)| 
&\ge&
|v(t)| - |B(t)-v(t)| \\
&\ge&
|B(t_0)| \cdot \left(\frac{t}{t_0}\right)^{-\mu} \cdot 
\left\{
\frac{1}{\sqrt{2}} - |\lambda |\cdot t^{1/2}_0\right\} 
\mbox{\hspace{2.55cm}by (\ref{3di21.3})}\\
&\ge&
\frac{1}{2} \cdot |B(t_0)| \cdot \left(\frac{t}{t_0}\right)^{-\mu}.
\end{eqnarray*}
Thus
\begin{eqnarray}
\int\limits^{t_2}_{t_{-2}} |B(t)|^2 dt 
&\ge&
\int\limits^{t_0}_{t_{-2}} |B(t)|^2 dt \nonumber \\
&\ge&
\frac{1}{4} \cdot |B(t_0)|^2 \cdot \frac{1}{2\mu-1} \cdot 
\left\{
\left(\frac{t_{-2}}{t_0}\right)^{-2\mu} \cdot t_{-2} - t_0\right\}
\nonumber \\
&\ge&
\frac{|B(t_0)|^2}{8(2\mu-1)} \cdot 
\left(\frac{t_0}{t_{-2}}\right)^{2\mu} \cdot t_{-2}
\label{g3.7}
\end{eqnarray}
because $t_0\ge t_{-1}\ge 2t_{-2}$ and thus 
$\left(\frac{t_0}{t_{-2}}\right)^{2\mu} - \left(\frac{t_0}{t_{-2}}\right) 
\ge \frac{1}{2} \left(\frac{t_0}{t_{-2}}\right)^{2\mu}$.

Hence by (\ref{g3.5}) and (\ref{g3.7})
$$
\frac{ 
\int\limits^{t_1}_{t_{-1}} |B(t)|^2 dt
      }
     {
\int\limits^{t_2}_{t_{-2}} |B(t)|^2 dt
      }
\le
\frac{\frac{8}{2\mu-1}\cdot |B(t_0)|^2 \cdot  
\left(\frac{t_1}{t_{-1}}\right)^{\mu} \cdot t_{1}}
{\frac{|B(t_0)|^2}{8(2\mu-1)} \cdot \left(\frac{t_0}{t_{-2}}\right)^{2\mu} 
\cdot t_{-2}} 
= 
2^6 \cdot \left(\frac{t_{-2}}{t_{-1}}\right)^{2\mu} \cdot \frac{t_1}{t_{-2}}.
$$
This finishes the proof of Proposition 3.2. \hfill$\Box$

\vspace{1cm}
\noindent
To conclude this section we cite the following proposition which has also 
been used in the proof of the Gluing Theorem.

\vspace{1cm}
{\bf Proposition 3.3.} {\it Let} $n\ge 3$. 
{\it Let} $a, b, c, \lambda\in\RR$ {\it such that} $a<b$ {\it and} $c>0$. 
{\it Let} $\rho:[a-c, b+c]\rightarrow\RR$ {\it be a smooth positive 
function such that}
$$
|\lambda| \cdot \|\rho\|_{L^{\infty}(a-c, b+c)} + \frac{1}{2} \cdot
\|\dot{\rho}\|_{L^{\infty}(a-c, b+c)} \le 1.
$$
{\it Then for every Dirac eigenspinor $\sigma$ for the eigenvalue}
$\lambda$ {\it on} $Z_{a-c, b+c} :=[a-c, b+c] \times S^{n-1}$ {\it with the 
warped product metric} $ds^2 = dt^2 + \rho(t)^2 d\sigma^2$ {\it the 
following estimate holds}
$$
\|\sigma\|^2_{L^2 (Z_{a, b})} \le \frac{b-a}{2c} 
\left\{
\|\sigma\|^2_{L^2 (Z_{b, b+c})} + \|\sigma\|^2_{L^2 (Z_{a-c, a})}
\right\}.
$$ 

\vspace{1cm}
\noindent
The proof can be found in \cite[Prop. 5.1]{b4}.

\vspace{2cm}

\section{Harmonic Sections}

In this last section we apply the Gluing Theorem to prove

\vspace{1cm}
{\bf Theorem 4.1.} {\it Let} $M$ {\it be a closed Riemannian manifold of 
dimension $n\equiv$ 3 mod 4.} 
{\it Let $D$ be an elliptic self-adjoint differential operator over $M$ 
of order 1.} 
{\it Let $U\subset M$ be a non-empty open subset}. 
{\it Let the restriction of $D$ to $U$ be a twisted Dirac operator}.

{\it Then one can deform the Riemannian metric in $U$ such that the 
resulting operator $\tilde{D}$ has non-trivial kernel}.

\vspace{1cm}
\noindent
Note that $\tilde{D}$ coincides with $D$ outside $U$ and that the 
connection of the coefficient bundle over $U$ is {\bf not} modified. 
Theorem 4.1 applies for example to the Dirac operators of spin,
$\Spinc$, or $\Spinh$ manifolds. 
By a suitable local deformation of the Riemannian metric while keeping 
the connection of the canonical line bundle fixed one can produce a 
non-trivial kernel for the Dirac operator of a $\Spinc$ manifold. 

Theorem 4.1 does not apply to the Euler operator $d + \delta$ acting
on forms even though $d + \delta$ is also a twisted Dirac operator 
(twisted by spinors).
The point is that in this case a change of the Riemannian metric also 
changes the connection on the coefficient bundle.

>From the discussion in the introduction we know that Theorem 4.1 is 
not true in dimension 2.
The classical Dirac operator on $S^2$ has no harmonic spinors no matter
which Riemannian metric on $S^2$ has been chosen.

\vspace{1cm}
{\bf Proof of Theorem 4.1.} Over $U$ our operator $D$ is a twisted Dirac 
operator, i.e. $D|_U = D^E$ where $E$ is some bundle over $U$ with 
connection $\nabla^E$. 
Choose a small $n$-ball $U_1 \subset\subset U$. 
By blowing up the metric $g$ in a neighborhood of $\overline{U}_1$ to 
some metric $\overline{g}$ one can make
$$
f=\mbox{ id} : (U_1, \overline{g}) \rightarrow (U_1, g)
$$
$\epsilon$-contracting for an arbitrarily small prescribed $\epsilon>0$.

\begin{center}
\begin{pspicture}(-6,-6)(6,13)

\psellipse[fillstyle=solid,fillcolor=lightgray](-5,0)(4,3)
\psellipse[fillstyle=solid,fillcolor=darkgray](-5,0)(2,1.5)
\rput(-5,0){\psframebox*[framearc=0.5]{$U_1$}}
\rput(-7.9,0){\psframebox*[framearc=0.5]{$U$}}
\rput(-5,-4){metric $g$}

\psellipse[fillstyle=solid,fillcolor=lightgray](5,0)(4,3)
\pscircle[fillstyle=solid,fillcolor=lightgray](5,6){6}
\pscustom[fillstyle=solid,fillcolor=darkgray]{
  \psarc*(5,6){6}{0}{180}
  \pscurve[liftpen=1](-1,6)(3,4)(7,4)(11,6)}
\pscustom[fillstyle=solid,fillcolor=lightgray]{
  \pscurve(3,0.4)(4,-0.2)(3.5,-0.6)(2.5,-0.8)
  \pscurve[liftpen=1](7.5,-0.8)(6.5,-0.6)(6,-0.2)(7,0.4)}
\psframe*[linecolor=lightgray](2.4,-0.9)(7.6,-0.7)
\rput(5,8){\psframebox*[framearc=0.5]{$U_1$}}
\rput(5,-1.5){\psframebox*[framearc=0.5]{$U$}}
\rput(5,-4){metric $\overline{g}$} 

\end{pspicture}
\begin{center}
{\bf Fig. 4}
\end{center}
\end{center} 

\vspace{1cm}
\noindent  
After trivializing the bundle $E$ over $U_1$ we can write
$$
\nabla^E = \partial + \Gamma
$$
where $\partial$ is a flat coordinate derivative and $\Gamma$ is an 
$\mbox{End}(E)$-valued $1$-form given by the Christoffel symbols.

For the pull-back of a connection by a map $f$ one has
$$
f^* \nabla^E = \partial + \Gamma \circ df.
$$
Since in our case $f=\mbox{ id}$ is $\epsilon$-contracting we see that 
with respect to the new metric $\overline{g}$ we have
$$
\|\Gamma\|_{L^{\infty}(U_1)} < \epsilon'
$$
for small $\epsilon'>0$. Choose a smaller $n$-ball $U_2 \subset\subset 
U_1$ and a cut-off function $\chi : M\rightarrow\RR$ such that
\begin{itemize}
\item[a)]
$0\le\chi\le 1$
\item[b)]
$\chi\equiv 0$ on $U_2$
\item[c)]
$\chi\equiv 1$ on $M-U_1$.
\end{itemize}
Define a new connection $\overline{\nabla}$ on $E$ by
$$
\overline{\nabla} := \partial+\chi\cdot\Gamma.
$$
Then $\overline{\nabla}$ is flat over $U_2$ and 
$\overline{\nabla}=\nabla^E$ over 
$U-U_1$. 
Denote by $D_1$ the twisted Dirac operator for the original connection 
$\nabla^E$ and the new metric $\overline{g}$, 
by $D_2$ the twisted Dirac operator for $\overline{\nabla}$ and 
$\overline{g}$.
We extend $D_1$ and $D_2$ to $M$ such that over $M-U$ we have 
$D_1=D_2=D$. Then the difference $D_1 - D_2$ is an operator of order $0$ 
with support in $U_1$. 
Its $L^2$-operator norm can be estimated by
$$
\|D_1 - D_2\|_{L^2, L^2} 
\le 
n \cdot \|\overline{\nabla} - \nabla^E\|_{L^{\infty}} 
<
n\cdot\epsilon' =: \epsilon''.
$$
In particular, $D_1$ and $D_2$ are $(\infty, \epsilon'')$-spectral close 
and $D_2$ is of Dirac type over $U_2$.

In \cite[Section 3]{b4} it is shown that in dimension $n\equiv 3 \mbox{ mod } 
4$ there exists a one-parameter family $g_{T}$ of Riemannian metrics on 
$S^n$, $T\in [a, b]$, such that the following holds for the associated 
Dirac operator $D_T$:
\begin{itemize}
\item[a)]
There is $\lambda(T)\in \mbox{ spec}\,(D_T)$ with $\lambda(a)=-1$, 
$\lambda(b)=+1$. 
\item[b)]
$\lambda(T)$ depends smoothly (actually linearly) on $T$.
\item[c)]
The multiplicity of $\lambda(T)$ is constant in $T$ and can be chosen 
arbitrarily large.
\item[d)]
$\lambda(T)$ is the only eigenvalue of $D_T$ in the interval $[-1, 1]$.
\end{itemize}
We choose the multiplicity of $\lambda(T)$ larger than the total number 
of eigenvalues of $D_2$ in the interval $[-1, 1]$. 
Now we can apply the Gluing Theorem to $D_2$ and $D_T$. We obtain a 
metric $\overline{g}_T$ on $M\# S^n=M$ such that the corresponding 
operator $\overline{D}_T$ on $M$ is $(1, \epsilon'')$-spectral close to 
$D_2\dot{\cup}D_T$. 

By construction of $\overline{g}_T$ the identity mapping 
$$
\mbox{ id} : (U_1, \overline{g}_T) \rightarrow (U_1, g)
$$
is still $\epsilon$-contracting. 
Thus $\tilde{D}_T$ and $\overline{D}_T$ are 
$(\infty, \epsilon'')$-spectral close where $\tilde{D}_T$ coincides with $D$ 
outside $U$ while over $U$, $\tilde{D}_T$ is the twisted Dirac operator 
for $\nabla^E$ and $\overline{g}_T$. 
It follows that $\tilde{D}_T$ is $(1, 2 \epsilon'')$-spectral close to 
$D_2\dot{\cup}D_T$. 

Since the multiplicity of $\lambda(T)$ was chosen larger than the total 
number of eigenvalues of $D_2$ in $[-1, 1]$, there are more negative 
eigenvalues of $D_2\dot{\cup}D_T$ in $[-1, 1]$ for $T=a$ whereas there 
are more positive ones for $T=b$. 
If $\epsilon''$ was chosen smaller than $\frac{1}{2}$, then this will 
remain true for the operator $\tilde{D}_T$. 
In particular, there must be some $T_0\in(a, b)$ for which $0$ is an 
eigenvalue. Hence Theorem 4.1 holds with the metric $\overline{g}_{T_0}$ 
and the corresponding operator $\tilde{D}=\tilde{D}_{T_0}$.
\hfill$\Box$

\vspace{2cm}

\vsp{1cm}

\noindent
Author's address:

\vsp{0.3cm}

\noindent
Mathematisches Institut\\
Universit\"at Freiburg\\
Eckerstr. 1\\
79104 Freiburg\\
Germany

\vsp{0.3cm}

\noindent
e-mail: {\typewriter baer@mathematik.uni-freiburg.de}

\end{document}